\documentclass[twoside]{article} 
\usepackage{graphicx}
\date{} 
\title{Asymptotic expansion of the Wright function for large variable and parameter}
\author{\sc R. B.\ Paris \\
{\em Division of Computing and Mathematics,} \\
{\em Abertay University, Dundee DD1 1HG, UK}}
\begin{document}
\def\f#1#2{\mbox{${\textstyle \frac{#1}{#2}}$}}
\def\dfrac#1#2{\displaystyle{\frac{#1}{#2}}}
\def\boldal{\mbox{\boldmath $\alpha$}}
\newcommand{\bee}{\begin{equation}}
\newcommand{\ee}{\end{equation}}
\newcommand{\sa}{\sigma}
\newcommand{\ka}{\kappa}
\newcommand{\al}{\alpha}
\newcommand{\la}{\lambda}
\newcommand{\ga}{\gamma}
\newcommand{\eps}{\epsilon}
\newcommand{\om}{\omega}
\newcommand{\fr}{\frac{1}{2}}
\newcommand{\fs}{\f{1}{2}}
\newcommand{\g}{\Gamma}
\newcommand{\br}{\biggr}
\newcommand{\bl}{\biggl}
\newcommand{\ra}{\rightarrow}
\newcommand{\gtwid}{\raisebox{-.8ex}{\mbox{$\stackrel{\textstyle >}{\sim}$}}}
\newcommand{\ltwid}{\raisebox{-.8ex}{\mbox{$\stackrel{\textstyle <}{\sim}$}}}
\renewcommand{\topfraction}{0.9}
\renewcommand{\bottomfraction}{0.9}
\renewcommand{\textfraction}{0.05}
\newcommand{\mcol}{\multicolumn}
\date{}
\maketitle
\pagestyle{myheadings}
\markboth{\hfill \sc R. B.\ Paris  \hfill}
{\hfill \sc Wright function with a large parameter\hfill}
\begin{abstract}
We consider the asymptotic expansion of the Wright function 
\[W_{\lambda,\mu}(z)=\sum_{n=0}^\infty\frac{z^n}{n! \g(\lambda n+\mu)}\qquad (\la>-1)\]
for large (positive and negative) variable and large parameter $\mu$. The analysis is based on use of the method of steepest descents applied to a suitable integral representation and, in part, complements the recent work of Ansari and Askari. 
Numerical results are presented to illustrate the accuracy of the different expansions obtained.
\vspace{0.3cm}

\vspace{0.1cm}
 
\end{abstract}

\vspace{0.3cm}

\noindent $\,$\hrulefill $\,$

\vspace{0.3cm}

\begin{center}
{\bf 1.\ Introduction}
\end{center}
\setcounter{section}{1}
\setcounter{equation}{0}
\renewcommand{\theequation}{\arabic{section}.\arabic{equation}}
The Wright function under consideration (also known as a generalised Bessel function) is defined by
\bee\label{e11}
W_{\lambda,\mu}(z)=\sum_{n=0}^\infty\frac{z^n}{n! \g(\lambda n+\mu)},
\ee
where $\lambda$ is supposed real and $\mu$ is, in general, an arbitrary complex parameter. The series 
converges for all finite $z$ provided $\lambda>-1$ and, when $\lambda=1$, it reduces to the modified Bessel function $z^{(1-\mu)/2}I_{\mu-1}(2\sqrt{z})$. This function has found recent application in the theory of fractional calculus \cite{GLM, L}.
The case corresponding to $\lambda=-\sigma$, $0<\sigma<1$ arises in the analysis of time-fractional diffusion and diffusion-wave equations \cite{MC}; see also \cite{PCM}. 
The case $\mu=0$ in (\ref{e11}) also finds application in probability theory and is discussed extensively in \cite{PV}, where it is denoted by 
$\phi(\lambda,0;z)=W_{\lambda,0}(z)$
and referred to therein as a  `reduced' Wright function.

The asymptotics of this function were first studied by Wright \cite{W34, W40} using the method of steepest descents applied to the integral representation
\bee\label{e12}
W_{\lambda,\mu}(z)=\frac{1}{2\pi i}\int_{-\infty}^{(0+)}t^{-\mu} \exp\,[t+zt^{-\lambda}]\,dt\qquad(\lambda>-1,\ \mu\in {\bf C});
\ee
see also the summaries in \cite{P, PCM} using a different approach. A recent paper by Ansari and Askari \cite{AA} has investigated the asymptotic expansion of (\ref{e11}) when the argument $z<0$ and the parameter $\mu>0$ are both large. More specifically, these authors considered the Wright function
\bee\label{e13}
{\cal W}^-_{\la,\nu}(x):=(\fs x)^\nu W_{\la,\nu+1}(-(\fs x)^{\la+1})=(\fs x)^\nu \sum_{n=0}^\infty \frac{(-)^n ( x/2)^{(\la+1)n}}{n! \g(\la n+\nu+1)}
\ee
for $x\to+\infty$ when $\nu=ax$, where $a>0$ is fixed, by employing the method of steepest descents applied to a suitable Laplace-type integral representation. Here we consider the same problem in more detail using the same approach. In a certain domain of the integration variable it is found that there are two relevant saddle points, which can either be real or a complex conjugate pair. We also discuss a special case when the parameters $a$ and $\la$ are connected in a specific manner that corresponds to the formation of a double saddle point.

The function of positive argument
\bee\label{e14}
{\cal W}^+_{\la,\nu}(x):=(\fs x)^\nu W_{\la,\nu+1}((\fs x)^{\la+1})=(\fs x)^\nu \sum_{n=0}^\infty \frac{( x/2)^{(\la+1)n}}{n! \g(\la n+\nu+1)}
\ee
is also considered for $x\to+\infty$ when $\nu=ax$, where $a>0$ is fixed. This function is described by a similar Laplace-type integral which, in contrast to that defining ${\cal W}^-_{\la,\nu}(x)$, can involve a more elaborate saddle-point structure. There is always a single real saddle but, dependent on the values of $a$ and $\la$, there can be several contributory complex saddles. In addition, it is possible to encounter
a Stokes phenomenon as the parameters $a$ and $\la$ are varied (with $x>0$). However, the additional complex saddle points yield contributions that are subdominant relative to that from the real saddle; these contributions may be neglected in certain applications.
\vspace{0.6cm}

\begin{center}
{\bf 2.\ The integral representation of ${\cal W}^-_{\la,\nu}(x)$}
\end{center}
\setcounter{section}{2}
\setcounter{equation}{0}
\renewcommand{\theequation}{\arabic{section}.\arabic{equation}}
From (\ref{e12}), we find the integral representation
\begin{eqnarray*}
{\cal W}^-_{\la,\nu}(x)&=&\frac{(\fs x)^\nu}{2\pi i}\int_{-\infty}^{(0+)} t^{-\nu-1} \exp\,[t-(x/2)^{\la+1}t^{-\la}]\,dt\\
&=&\frac{1}{2\pi i}\int_{-\infty}^{(0+)} \tau^{-\nu-1} \exp\,[(x/2) (\tau-\tau^{-\la})]\,d\tau,
\end{eqnarray*}
where $t=\fs x\tau$. If we make the further change of variable $\tau=e^u$ and put $\nu=ax$, we obtain
\bee\label{e21}
{\cal W}^-_{\la,\nu}(x)=\frac{1}{2\pi i}\int_C e^{xh(u)} du,\qquad h(u):=\fs(e^u-e^{-\la u})-au.
\ee
This last transformation causes all the Riemann sheets in the $\tau$-plane to appear as horizontal strips of width $2\pi$ in the $u$-plane.
The loop in the $\tau$-plane can be taken to be a circle of radius $\rho$ about the origin together with straight line segments on the upper and lower sides of the branch cut on the negative real axis. The map $C$ of this path in the $u$-plane then consists of the three sides of the rectangle with vertices at $\infty-\pi i$, $\log\,\rho-\pi i$, $\log\,\rho+\pi i$ and $\infty+\pi i$; see Fig.~1. An equivalent version of (\ref{e21}) was given in \cite{AA}.
\begin{figure}[th]
	\begin{center}{($a$)}\ \includegraphics[width=0.30\textwidth]{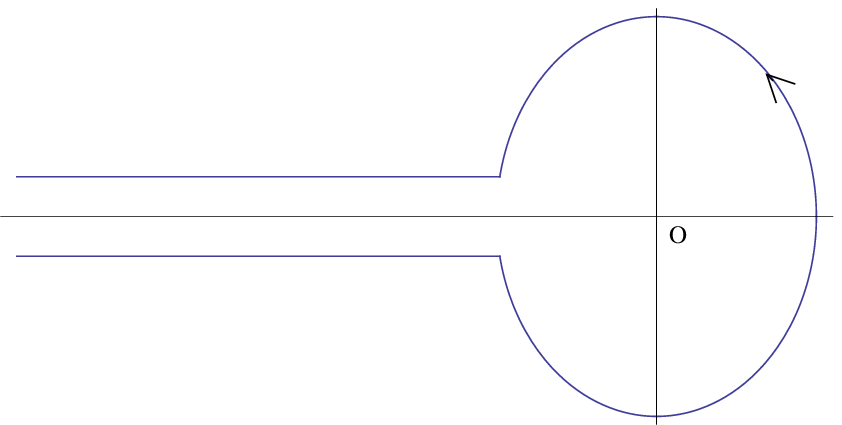}\qquad\qquad{($b$)}\includegraphics[width=0.30\textwidth]{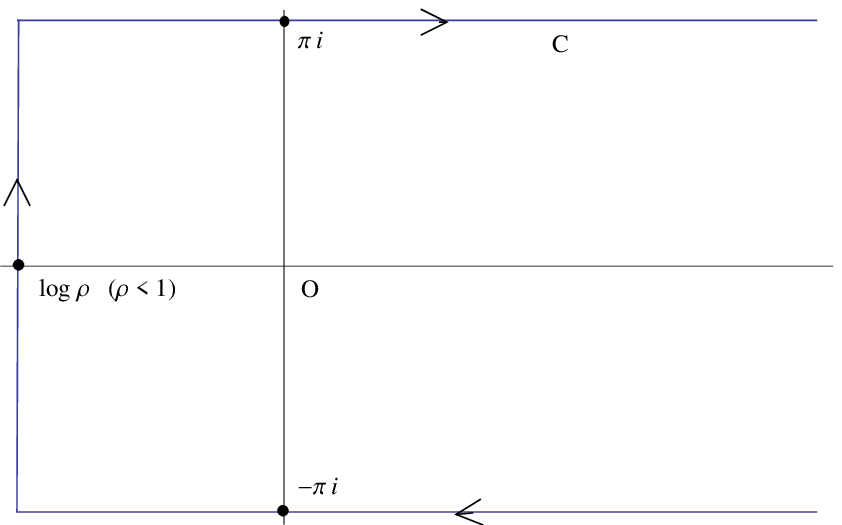}
\caption{\small{($a$) The loop in the $\tau$-plane with a circular path round the origin of radius $\rho$; ($b$) the contour $C$ in the $u$-plane shown with $\rho<1$.}}
\end{center}
\end{figure}

Saddle points of the integrand in (\ref{e21}) occur at $h'(u)=0$; that is, at the point $u=u_0$ of the equation
\bee\label{e22}
e^{u_0}+\la e^{-\la u_0}=2a.
\ee
It is sufficient to confine our attention to the strip $-\pi\leq \Im (u) \leq\pi$ when considering the location of the saddle points.
By inspection of the above equation, it is easily seen that when $-1<\la<0$ there is only one real saddle in this strip
and when $\la>0$, there are two saddles, either both real or a complex conjugate pair. The saddles coincide to form a double saddle point when $h'(u)=h''(u)=0$ at
\bee\label{e23}
u_0=\frac{2\log\,\la}{1+\la}.
\ee
From (\ref{e22}), this requires the relation between the parameters $a$ and $\la$ given by
\bee\label{e23a}
a=\frac{(1+\la)}{2}\,\la^\gamma,\qquad \gamma:=\frac{1-\la}{1+\la}.
\ee
A plot of this curve is shown in Fig.~2. Above this curve the saddles are real and below they are a complex conjugate pair; on the curve there is a double saddle. Routine calculations show that the maximum of the curve occurs at $a\doteq 1.19123$ when $\la\doteq 2.09350$.

\begin{figure}[th]
	\begin{center}	\ \includegraphics[width=0.50\textwidth]{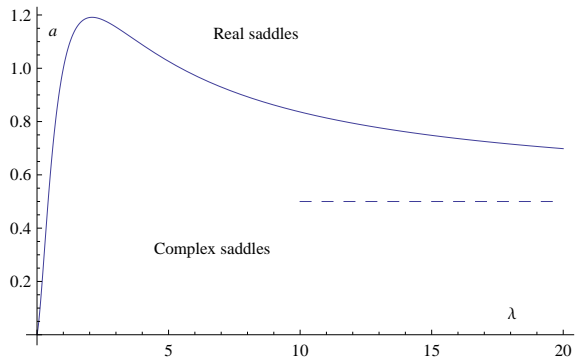}
\caption{\small{Plot of the curve in (\ref{e23a}) corresponding to the formation of a double saddle on the real $u$-axis when $\la>0$. Above this curve the two saddles are real and below they form a complex conjugate pair. The dashed line $a=\frac{1}{2}$ represents the asymptotic limit of the curve. }}
\end{center}
\end{figure}

Typical paths of steepest descent through the relevant saddles in the strip $|\Im (u)|\leq\pi$ are shown in Fig.~3. In each case the paths that 
pass through the contributory saddle(s) pass to infinity in the right half-plane with $\Im (u)=\pm\pi$. When $\la>0$ and $(a, \la)$ is situated above the curve in Fig.~2 only the larger real saddle contributes to the expansion of
${\cal W}^-_{\la,\nu}(x)$; the path through the smaller saddle is a path of steepest ascent. The integration path $C$ can, in each case, be deformed to pass over the steepest descent paths passing through the contributory saddle(s); in the case of complex saddles this can be achieved by letting $\rho\to0$. 
\begin{figure}[th]
	\begin{center}{($a$)}\ \includegraphics[width=0.30\textwidth]{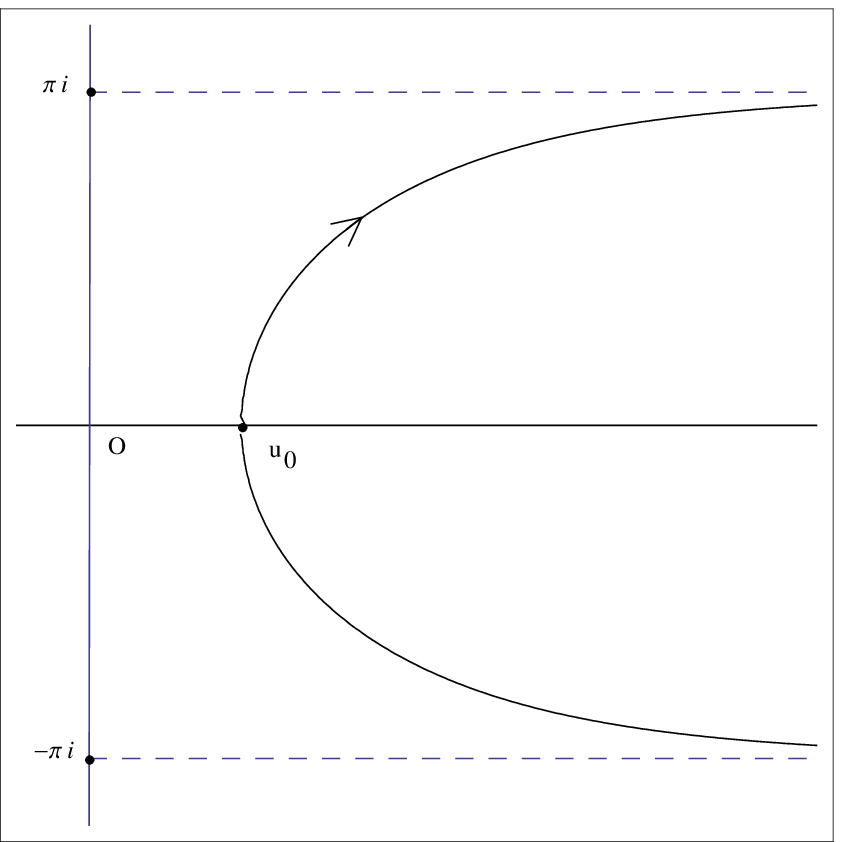}\qquad\qquad{($b$)}\ \includegraphics[width=0.30\textwidth]{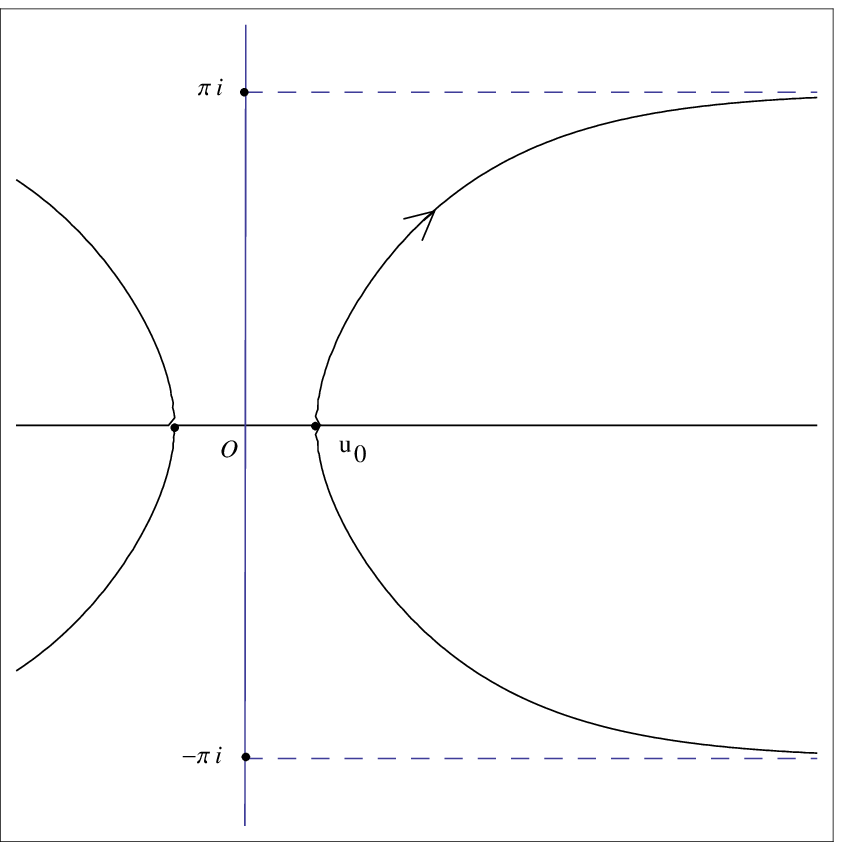}\\
\vspace{0.3cm}
	
	($c$)\ \includegraphics[width=0.30\textwidth]{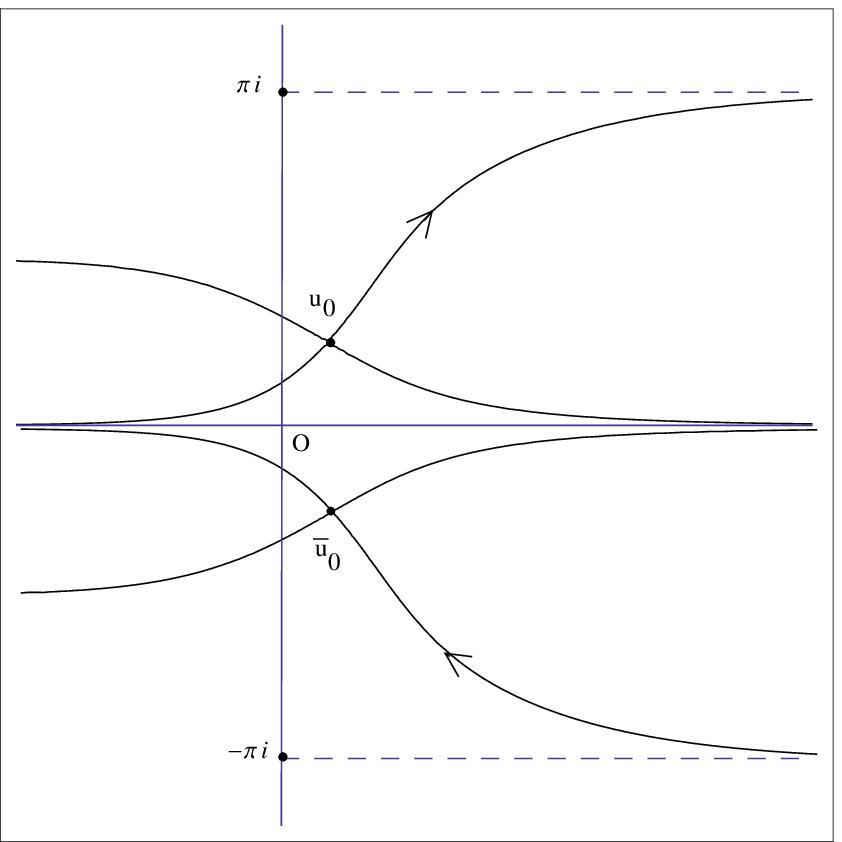}\qquad\qquad{($d$)}\ \includegraphics[width=0.30\textwidth]{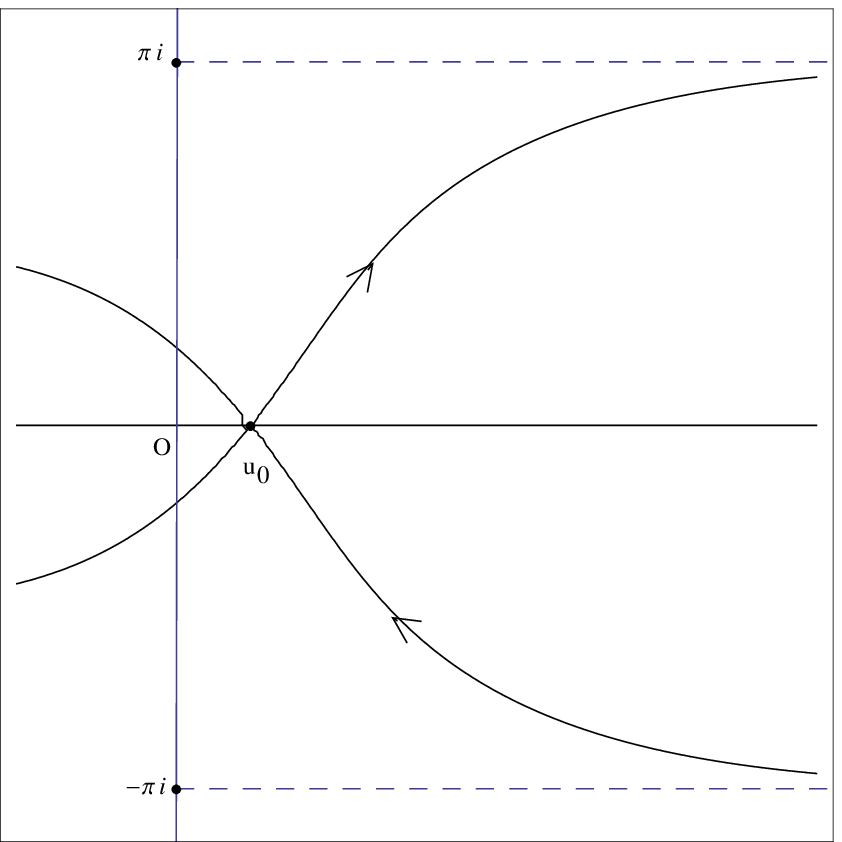}
\caption{\small{ Steepest descent paths ($a$) when $-1<\la<0$, ($b$) real saddle case when $\la>0$; ($c$) complex saddles case and ($d$) double saddle case. The arrows denote the direction of integration.}}
\end{center}
\end{figure}

\vspace{0.6cm}
\begin{center}
{\bf 3.\ Asymptotic expansions}
\end{center}
\setcounter{section}{3}
\setcounter{equation}{0}
\renewcommand{\theequation}{\arabic{section}.\arabic{equation}}
In this section we consider the asymptotic expansions of ${\cal W}^-_{\la,\nu}(x)$ in the three cases that can arise, namely a real contributory saddle, a complex conjugate pair of saddles and a real double saddle.
\vspace{0.3cm}

\noindent{\bf 3.1\ \ Real saddle point}
\vspace{0.2cm}

\noindent
When $-1<\la<0$ and $\la>0$ with $(a,\la)$ situated above the curve in Fig.~2, the contributory saddle is real; see Fig.~3(a, b). We have from (\ref{e21})
\[{\cal W}^-_{\la,\nu}(x)=\frac{e^{xh(u_0)}}{2\pi i} \int_{-\infty}^\infty e^{-xw^2/2}\,\frac{du}{dw}\,dw,\]
where we have introduced the new variable $w$ by
\[-\frac{1}{2}w^2=h(u)-h(u_0)=\frac{h_0''}{2!} (u-u_0)^2+\frac{h_0'''}{3!} (u-u_0)^3+\cdots\]
and for brevity we denote the derivatives of $h(u)$ evaluated at $u=u_0$ by $h_0^{(n)}\equiv h^{(n)}(u_0)$.
Inversion of this expansion using the {\tt InverseSeries} command in {\it Mathematica} (essentially Lagrange inversion) yields
\bee\label{e300}
u-u_0=\frac{w}{\sqrt{-h_0''}}+\frac{h_0'''w^2}{6h_0''^2}-\frac{(5h_0'''^2-3h_0'' h_0^{iv})}{72 h_0''^3 \sqrt{-h_0''}}\,w^3+\cdots\ .
\ee
Upon differentiation this yields
\[\frac{du}{dw} \stackrel{e}{=} \frac{1}{\sqrt{-h_0''}} \sum_{k\geq0} (-)^k A_k w^{2k},\]
where $\stackrel{e}{=}$ denotes the inclusion of {\em only the even powers of $w$}, since odd powers will not enter into this calculation. We have 
\bee\label{e2h}
h_0''=\fs(e^{u_0}-\la^2 e^{-\la u_0})=\fs(1+\la)e^{u_0}-a\la
\ee
and note that $0<h_0''<a$. The first few coefficients $A_k$ are: 
\begin{eqnarray}
A_0&=&1, \qquad A_1=\frac{1}{24 h_0''}(5{\bf h}_3^2-3{\bf h}_4),\nonumber\\
A_2&=&\frac{1}{3456h_0''^2}(385{\bf h}_3^4-620{\bf h}_3^2 {\bf h}_4+105{\bf h}_4^2+168{\bf h}_3{\bf h}_5-24{\bf h}_6)\\
A_3&=&\frac{1}{6220800h_0''^3}(425425{\bf h}_3^6-1126125{\bf h}_3^4{\bf h}_4+675675{\bf h}_3^2{\bf h}_4^2-51975{\bf h}_4^3+360360{\bf h}_3^3{\bf h}_5\nonumber\\
&& -249480{\bf h}_3{\bf h}_4{\bf h}_5+13608{\bf h}_5^2-83160{\bf h}_3^2{\bf h}_6+22680{\bf h}_4{\bf h}_6+12960{\bf h}_3{\bf h}_7-1080{\bf h}_8)\nonumber,
\end{eqnarray}
where we have defined
\[{\bf h}_n:=\frac{h_0^{(n)}}{h_0''},\qquad h_0^{(n)}= \fs(1-(-\la)^{n-1})e^{u_0}+a(-\la)^{n-1}\qquad (n\geq 3).\]

Then
\[{\cal W}^-_{\la,\nu}(x)\sim \frac{e^{xh(u_0)}}{\pi\sqrt{2xh_0''}} \sum_{k\geq 0}\frac{(-)^k A_k}{(x/2)^k} \int_0^\infty e^{-s} s^{k-1/2} ds\]
\bee\label{e32}
= \frac{e^{xh(u_0)}}{\sqrt{2\pi xh_0''}}\sum_{k\geq0} \frac{(-)^k (\fs)_k A_k}{(x/2)^k}\hspace{0.3cm}
\ee
as $x\to+\infty$ with $\nu=ax$, where $(a)_n=a(a+1) \ldots (a+n-1)$ is the Pochhammer symbol. An expansion equivalent to this is given in Theorem 2.2 in \cite{AA}. 
\vspace{0.3cm}

\noindent{\bf 3.2\ \ Complex saddle points}
\vspace{0.2cm}

\noindent
When $\la>0$ and $(a, \la)$ is situated below the curve in Fig.~2 the contributory saddles are complex as illustrated in Fig.~3(c). The contribution to ${\cal W}^-_{\la,\nu}(x)$ from the upper steepest descent path is given by (compare (\ref{e32}))
\[\frac{e^{xh(u_0)}}{i\sqrt{2\pi xh_0''}}\sum_{k\geq0} \frac{(-)^k (\fs)_k A_k}{(x/2)^k},\]
where in this case $h_0''$ given in (\ref{e2h}) is complex-valued. Since the contribution
from the lower path yields the conjugate expression,  we therefore obtain
\bee\label{e33}
{\cal W}^-_{\la,\nu}(x) \sim \sqrt{\frac{2}{\pi x}} \,\Re \bl\{\frac{e^{xh(u_0)}}{\sqrt{h_0''}} \sum_{k\geq0} \frac{(-)^k (\fs)_k A_k}{(x/2)^k}\br\}
\ee
as $x\to+\infty$ with $\nu=ax$. The leading term of this expansion is given in an equivalent form in Theorem 3.1 in \cite{AA}.

It is worth remarking that in \S\S 3.1, 3.2 the exact location of the saddle $u_0$ is not available in closed form. This necessitates the numerical solution of (\ref{e22}) to obtain $u_0$ when the parameters $a$ and $\la$ have specific values. It is apparent that the presentation of higher coefficients $A_k$ becomes prohibitive on account of their rapidly increasing complexity. This point is discussed further in Section 5 where it is indicated how more coefficients $A_k$ can be generated by numerical reversion in specific cases.
\vspace{0.3cm}

\noindent{\bf 3.3\ \ Double saddle point}
\vspace{0.2cm}

\noindent
In the case of the double saddle point we are in the fortunate position of having an exact expression for the location of this saddle given in (\ref{e23}) combined with the condition (\ref{e23a}). The contribution from the upper half of the steepest descent path illustrated in Fig.~3(d) is
\[I=\int_{u_0}^{\infty+\pi i}e^{xh(u)}du=e^{xh(u_0)}\int_0^\infty e^{-xw^3/3} \frac{du}{dw}\,dw,\]
where we have introduced the new variable $w$ by
\[-\frac{w^3}{3}=h(u)-h(u_0)=\frac{h_0'''}{3!}(u-u_0)^3+ \frac{h_0^{(4)}}{4!} (u-u_0)^4+\cdots\ .\]
Inversion of this expansion (using {\it Mathematica}) yields for the upper half of the integration path
\[u-u_0=\frac{2^{2/3}e^{\pi i/3}}{H^{1/3}}\,w+\frac{(\la-1)e^{2\pi i/3}}{2^{2/3}\cdot 3 H^{2/3}}\,w^2-\frac{(1-6\la+\la^2)}{60H}\,w^3+\cdots\ ,\]
where $H:=2h_0'''$.
Upon differentiation, we then obtain 
\[\frac{du}{dw}=\frac{2^{2/3}}{H^{1/3}}\sum_{k\geq0} \frac{e^{\frac{1}{3}\pi i(k+1)}B_k}{H^{k/3}}\,w^k,\]
where
\begin{eqnarray*}
B_0&=&1,\quad B_1=\frac{\la-1}{2^{1/3}\cdot 3},\quad B_2=\frac{1-6\la+\la^2}{2^{2/3}\cdot 20},\quad B_3=\frac{1}{1620}(5+93\la-93\la^2-5\la^3),\\
B_4&=&-\frac{1}{2^{1/3}\cdot 136080}(277+826\la-6114\la^2+826\la^3+277\la^4),\\
B_5&=&\frac{1}{2^{2/3}\cdot 16800}(1-61\la-254\la^2+254\la^3+61\la^4-\la^5),\\
B_6&=&\frac{1}{10497600}(959+7098\la-2031\la^2-58708\la^3-2031\la^4+7098\la^5+959\la^6).
\end{eqnarray*}

The contribution from the upper half of the steepest descent path then becomes
\[I\sim \frac{2^{2/3} e^{xh(u_0)}}{3(Hx/3)^{1/3}} \sum_{k\geq0} \frac{e^{\frac{1}{3}\pi i(k+1)}B_k}{(Hx/3)^{k/3}} \int_0^\infty e^{-s} s^{k/2-2/3}\,ds\]
\[=\frac{2^{2/3} e^{xh(u_0)}}{3(Hx/3)^{1/3}} \sum_{k\geq0} \frac{e^{\frac{1}{3}\pi i(k+1)}B_k}{(Hx/3)^{k/3}}\,\g\bl(\frac{k+1}{3}\br).\]
The contribution from the lower half of the integration path is given by the conjugate expansion. Hence, when $\nu=ax$, $\la>0$ and the parameter $a$ satisfies the condition (\ref{e23a}), we obtain the expansion
\begin{eqnarray}
{\cal W}^-_{\la,\nu}(x)\!\!&\sim&\!\! \frac{1}{2\pi i}(I-{\overline I})\nonumber\\
&=&\!\!\frac{2^{2/3} e^{xh(u_0)}}{3\pi(Hx/3)^{1/3}} \sum_{k\geq0} \frac{B_k}{(Hx/3)^{k/3}}\,\g\bl(\frac{k+1}{3}\br) \sin \frac{\pi}{3} (k+1)\label{e31}
\end{eqnarray}
as $x\to+\infty$, where 
\[H=\la^{2/(1+\la)}(1+\la),\qquad h(u_0)=\frac{(\la^2-1)}{2\la^2} e^{u_0}-\la^\gamma \log\,\la,\]
with $u_0$ and $\gamma$ given in (\ref{e23}) and (\ref{e23a}).
We observe that the terms corresponding to $k=2, 5, 8, \ldots$ make no contribution to the expansion (\ref{e31}).

\vspace{0.6cm}

\begin{center}
{\bf 4.\ The expansion of ${\cal W}^+_{\la,\nu}(x)$}
\end{center}
\setcounter{section}{4}
\setcounter{equation}{0}
\renewcommand{\theequation}{\arabic{section}.\arabic{equation}}
We consider the function ${\cal W}^+_{\la,\nu}(x)$ defined in (\ref{e14}), which has the integral representation 
\bee\label{e41}
{\cal W}^+_{\la,\nu}(x)=\frac{1}{2\pi i}\int_C e^{x{\tilde h}(u)}du,\qquad {\tilde h}(u):=\fs(e^u+e^{-\la u})-au,
\ee
where $C$ is the same path as in (\ref{e21}). Paths of steepest descent (resp. ascent) as $\Re (u)\to-\infty$ asymptote to $\Im (u)=\pi n/\la$ for odd (resp. even) $n$. 
Saddle points are given by
\[e^{u_0}-\la e^{-\la u_0}=2a;\]
when $-1<\la<0$, there is just one real saddle $u_0$ in the principal strip $-\pi\leq\Im (u)\leq\pi$. When $\la>0$, there is the real saddle $u_0$ together with an infinite string of complex saddles $u_k=X_k\pm iY_k$, $k=1, 2, \ldots$ approximately parallel to the imaginary axis with $X_k>0$. The saddle $u_1$ has $\pi/\la<Y_1<3\pi/\la$ while it is found that $Y_k\simeq (2k-1)\pi/\la$ as $k$ increases. As $\la$ increases, more of these complex saddles enter the principal strip, some of which may, dependent on the values of $a$ and $\la$, interact with the steepest descent path through $u_0$. We define the second derivative ${\tilde h}_k''$ by
\[{\tilde h}_k'':=\fs(e^{u_k}+\la^2 e^{-\la u_k})=\fs(1+\la) e^{u_k}-a\la \qquad(k=0, 1, 2 \ldots).\]
\begin{figure}[th]
	\begin{center} \includegraphics[width=0.50\textwidth]{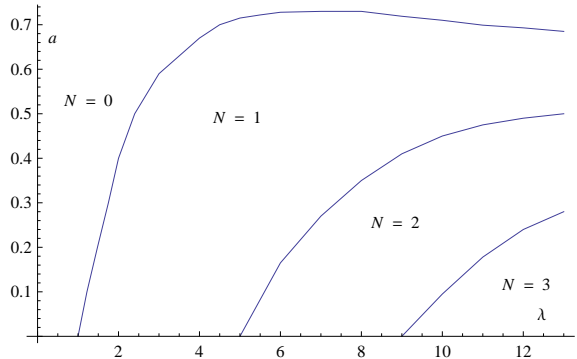}
\caption{Curves in the $a,\la$-plane on which an additional pair of complex saddles appears (via a Stokes phenomenon).}
\end{center}
\end{figure}

For $a$ and $\la$ in a certain domain (labelled $N=0$ in Fig.~4), and also for $-1<\la<0$, the steepest descent path passes through only the saddle $u_0$ on the real axis, as shown in Fig.~5(a). In this case the expansion of ${\cal W}^+_{\la,\nu}(x)$ is given by
\[{\cal W}^+_{\la,\nu}(x)\sim I_0(a,\la;x):=\frac{e^{x{\tilde h}(u_0)}}{\sqrt{2\pi x {\tilde h}_0''}} \sum_{k=0}^\infty \frac{(-)^k (\fs)_kB_k^{(0)}}{(x/2)^k}\qquad (x\to+\infty),\]
where the coefficients $B_k^{(0)}$ are obtained in the same manner as the $A_k$ described in Section 3, with $h(u)$ replaced by ${\tilde h}(u)$. As $\la$ ($>1$) increases in the domain labelled $N=1$ in Fig.~4, the steepest descent path through $u_0$ either passes to $\infty\pm\pi i$ (as shown in Fig.~5(a)) or to $-\infty$ along the directions $\Im (u)=\pm\pi/\la$ and thence over the saddles $u_1$, ${\overline u}_1$ to the endpoints $\infty\pm\pi i$; see Fig.~5(c).
The intermediate case shown in Fig.~5(b) shows the steepest descent path through $u_0$ connecting with the adjacent saddles $u_1$ and ${\overline u}_1$ to produce a Stokes phenomenon. As $\la$ increases further in a certain domain of the $a, \la$-plane, the saddles $u_2$, ${\overline u}_2$ become connected; see Fig. 5(d). In each case, the integration path $C$ can be deformed to coincide with these different steepest descent paths. It is worth noting that the appearance of each new pair of complex saddles results in a Stokes phenomenon. We do not consider the details of this transition here.
\begin{figure}[th]
	\begin{center}{($a$)}\ \includegraphics[width=0.30\textwidth]{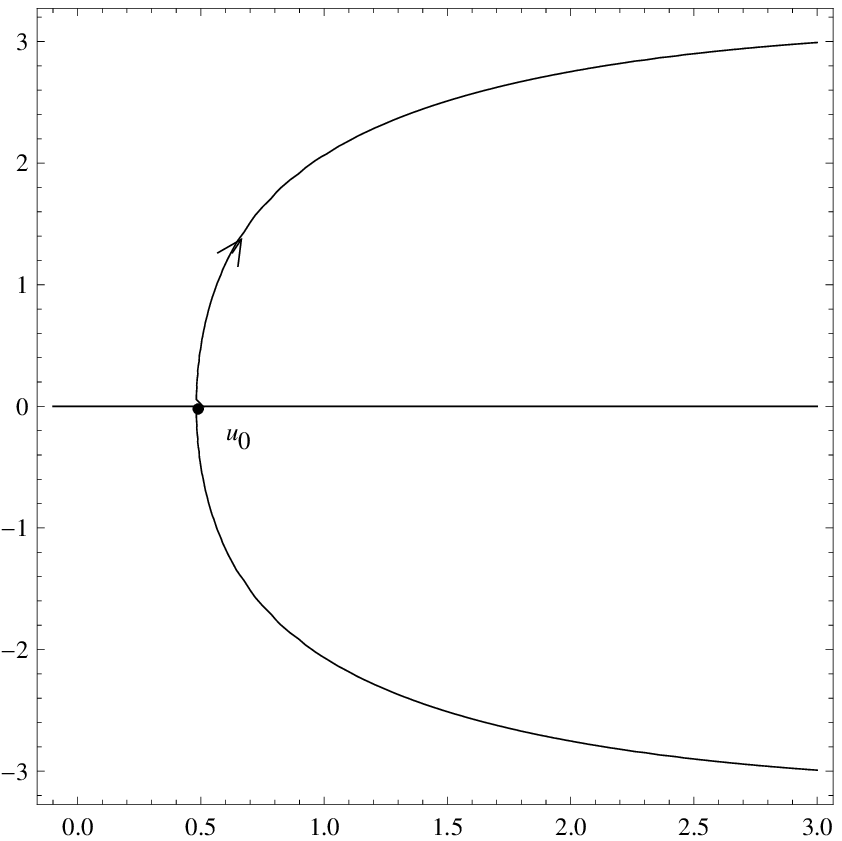}\qquad\qquad{($b$)}\ \includegraphics[width=0.30\textwidth]{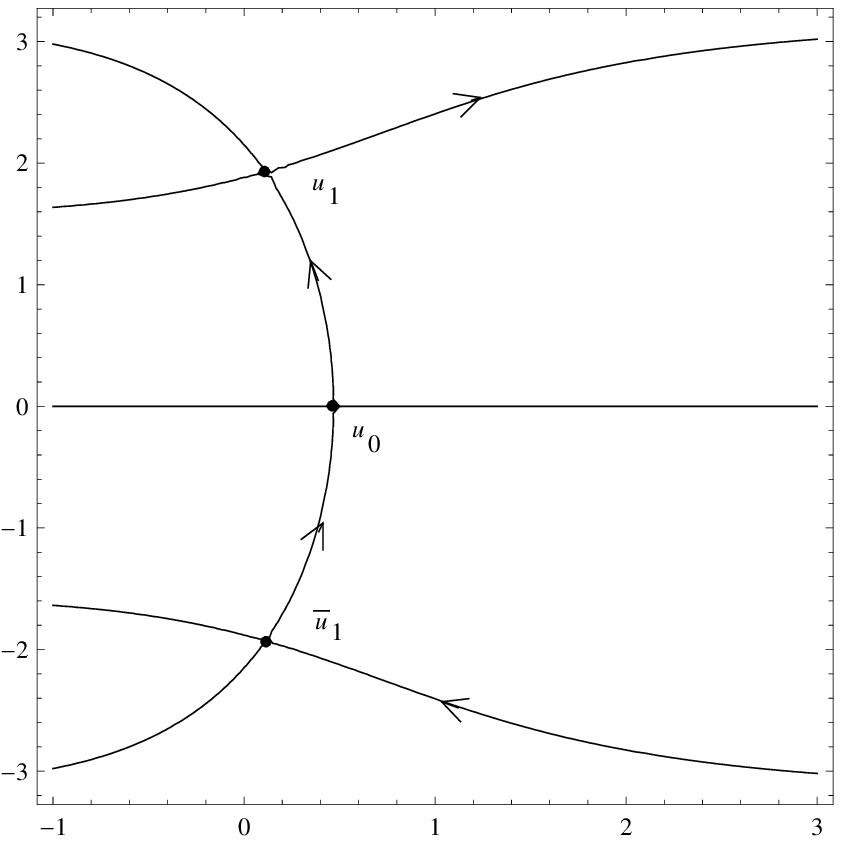}\\
\vspace{0.3cm}
	
	($c$)\ \includegraphics[width=0.30\textwidth]{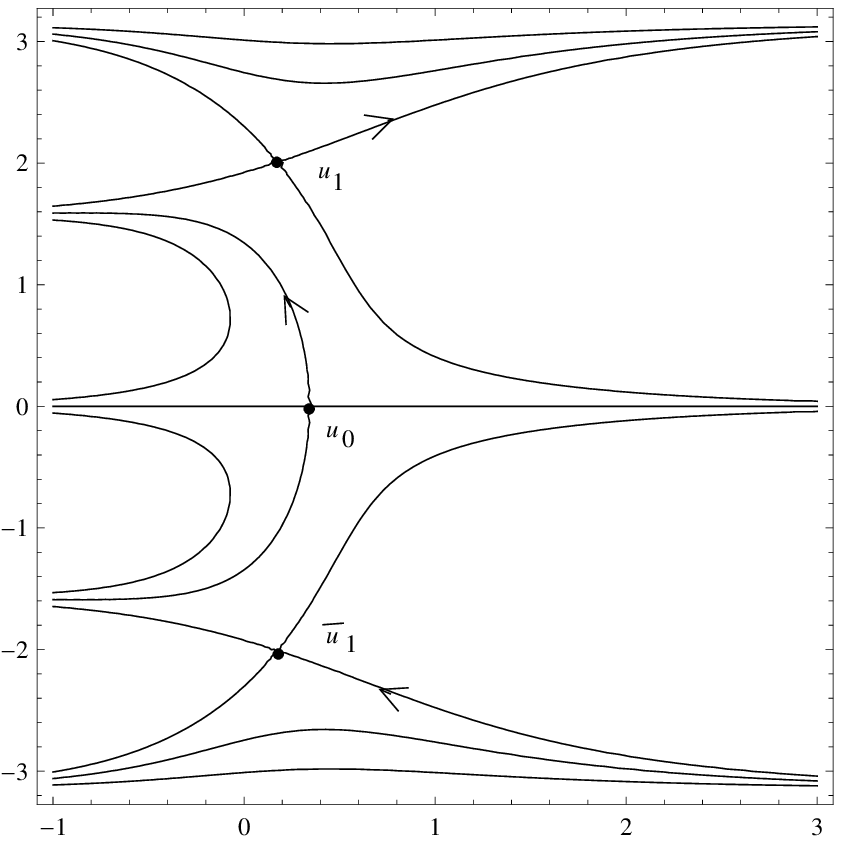}\qquad\qquad{($d$)}\ \includegraphics[width=0.30\textwidth]{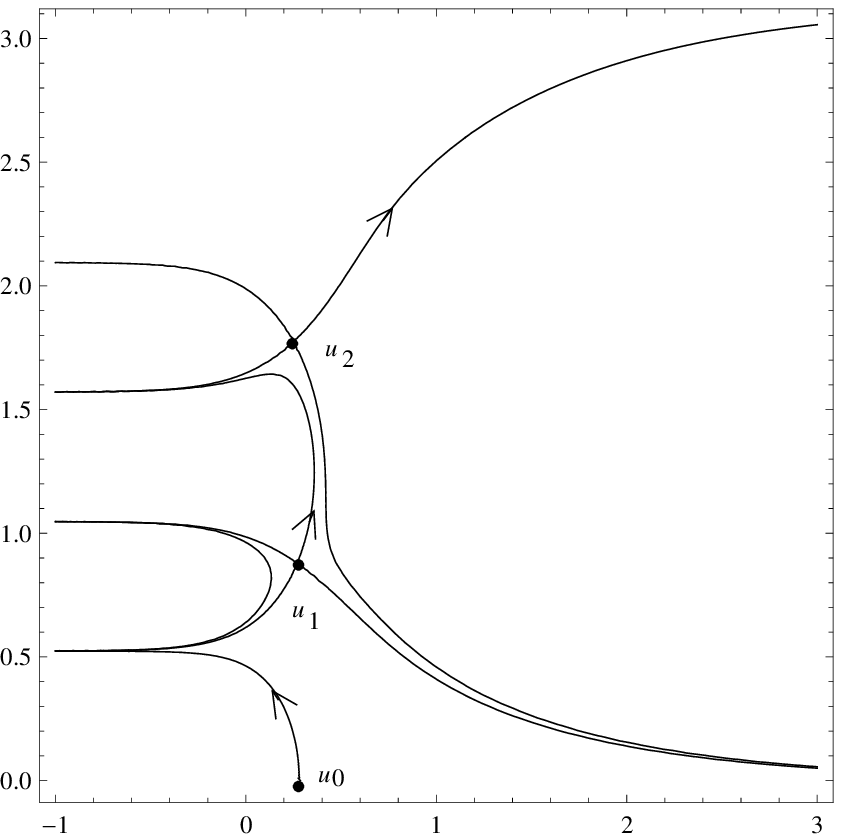}
\caption{\small{ Steepest descent paths ($a$) when $\la=1, a=0.50$ ($N=0$), ($b$) $\la=2, a=0.4075$ (on the boundary between $N=0$ and $N=1$ in Fig.~4); ($c$) $\la=2, a=0.60$ ($N=1$) and ($d$) $\la=6, a=0.10$ ($N=2$). For clarity the paths in ($d$) are shown only in the upper-half plane; a symmetrical distribution is present in the lower-half plane. The arrows indicate the direction of integration.}}
\end{center}
\end{figure}
\begin{table}[th]
\caption{\footnotesize{Values of the coefficients $A_k$ and the absolute relative error in ${\cal W}^-_{\la,\nu}(x)$ resulting from the expansion (\ref{e32}) for different truncation index $k$ when $x=40$. The value of the real saddle $u_0$ is indicated. }}
\begin{center}
\begin{tabular}{|l||l|l||r|r||r|r|}
\hline
\mcol{1}{|c||}{} & \mcol{2}{c||}{$\la=-0.25,\ a=1$} & \mcol{2}{c||}{$\la=1,\ a=1.20$} & \mcol{2}{c|}{$\la=0.50,\ a=0.80$}\\
\mcol{1}{|c||}{} & \mcol{2}{c||}{$u_0=0.83644438$}& \mcol{2}{c||}{$u_0=0.62236250$} & \mcol{2}{c|}{$u_0=0.12181472$}\\
\mcol{1}{|c||}{$k$} & \mcol{1}{c|}{$A_k$} & \mcol{1}{c||}{Error} & \mcol{1}{c|}{$A_k$} & \mcol{1}{c||}{Error} & \mcol{1}{c|}{$A_k$} & \mcol{1}{c|}{Error}\\
\hline
&&&&&&\\[-0.3cm]
0  & 1.000000 & $2.019(-3)$ & 1.000000 & $1.839(-2)$ & 1.000000 & $1.331(-2)$\\
1  & $+8.087175(-2)$& $3.189(-6)$ & 0.839435 & $2.655(-3)$ & 0.571373 & $9.888(-4)$\\
2  & $+1.681574(-3)$ & $2.995(-8)$ & 1.770726 & $7.334(-4)$ & 0.598231 & $1.490(-4)$\\
3  & $-1.284463(-4)$ & $2.168(-10)$& 4.345560 & $3.037(-4)$ & 0.768780 & $3.359(-5)$\\
4  & $-5.177287(-6)$ & $4.055(-12)$& 11.283213& $1.678(-4)$ & 1.050527 & $1.008(-5)$\\
5  & $+4.453244(-7)$ & $6.262(-14)$& 30.237515& $1.164(-4)$ & 1.483045 & $3.788(-6)$\\
\hline
\end{tabular}
\end{center}
\end{table}

The contribution from the pair of complex saddles $u_j$ and ${\overline u}_j$ is 
\[I_j(a, \la;x):=\sqrt{\frac{2}{\pi x}} \,\Re \bl\{\frac{e^{x{\tilde h}(u_j)}}{\sqrt{{\tilde h}_j''}} \sum_{k\geq0} \frac{(-)^k (\fs)_k B^{(j)}_k}{(x/2)^k}\br\}\qquad(x\to+\infty),\]
where the coefficients $B_k^{(j)}$ are computed as in (\ref{e33}).
Hence the expansion of ${\cal W}^+_{\la,\nu}(x)$ takes the form
\bee\label{e42}
{\cal W}^+_{\la,\nu}(x)\sim \sum_{j=0}^N I_j(a,\la;x)\qquad (x\to+\infty),
\ee
where $N$ denotes the number of contributory pairs of complex saddles. The different values of $N$ in the $a, \la$-plane are shown in Fig.~4, where each curve represents the appearance (via a Stokes phenomenon) of the pair of contributory complex saddles $u_N$ and ${\overline u}_N$. The dominant contribution results from $I_0(a,\la;x)$
with the series $I_j(a,\la;x)$ being progressively less dominant as $j$ increases. It is found numerically that the final series $I_N(a,\la;x)$ can be either exponentially large or small as $x\to+\infty$ according to the location of the parameters in the $a, \la$-plane.

\vspace{0.6cm}

\begin{center}
{\bf 5.\ Numerical verification}
\end{center}
\setcounter{section}{5}
\setcounter{equation}{0}
\renewcommand{\theequation}{\arabic{section}.\arabic{equation}}
In this section we present some numerical results that confirm the accuracy of the different expansions obtained for ${\cal W}^\pm_{\la,\nu}(x)$. We first present results for ${\cal W}^-_{\la,\nu}(x)$. In Table 1 we show the absolute relative error\footnote{In the tables we have adopted the convention of writing $x(y)$ for $x\times 10^y$.} in the computation of ${\cal W}^-_{\la,\nu}(x)$ using the expansion in (\ref{e32}) for different values of the parameters $a$ and $\la$ as a function of the truncation index $k$. The value of the real saddle point $u_0$ is given together with the first six coefficients $A_k$. The coefficients with $k\leq 3$ can be obtained from (\ref{e32}), but it is more direct, with specific values of $a$ and $\la$, to compute these coefficients using the {\tt InverseSeries} command in {\it Mathematica} to generate the numerical equivalent of (\ref{e300}).
\begin{table}[h]
\caption{\footnotesize{Values of the coefficients $A_k$ and the absolute relative error in ${\cal W}^-_{\la,\nu}(x)$ resulting from the expansion (\ref{e32}) for different truncation index $k$ when $x=40$, $\la=1.50$, $a=0.50$.}}
\begin{center}
\begin{tabular}{|l|l|r|}
\hline
\mcol{1}{|c|}{$k$} & \mcol{1}{c|}{$A_k$} & \mcol{1}{c|}{Error}\\
\hline
&&\\[-0.3cm]
0  & 1.000000                    & $6.233(-3)$ \\
1  & $+0.00929936 + 0.19815193i$ & $1.157(-4)$ \\
2  & $-0.08194718 + 0.01105633i$ & $1.416(-5)$ \\ 
3  & $-0.00729013 - 0.04233881i$ & $5.787(-7)$ \\
4  & $+0.02361754 - 0.00432441i$ & $1.840(-7)$ \\
5  & $+0.00253174 + 0.01363033i$ & $1.014(-8)$ \\
\hline
\end{tabular}
\end{center}
\end{table}

Table 2 shows the errors involved in using the expansion (\ref{e33}) when the contributory saddles are a conjugate pair located at $u_0=0.24834557 \pm 0.90919096i$. The coefficients $A_k$ in the case are complex valued. The special case when $a$ and $\la$ are linked via (\ref{e23a}) corresponding to a double saddle point is shown in Table 3 for different $\la$. All cases presented correspond to sub-optimal truncation as the relative error is seen to steadily decrease with increasing truncation index $k$.
\begin{table}[h]
\caption{\footnotesize{Values of the absolute relative error in ${\cal W}^-_{\la,\nu}(x)$ in the double saddle case resulting from the expansion (\ref{e31}) for different $\la$ and truncation index $k$ when $x=40$.}}
\begin{center}
\begin{tabular}{|l|lll|}
\hline
&&&\\[-0.3cm]
\mcol{1}{|c|}{$k$} & \mcol{1}{c}{$\la=0.50$} & \mcol{1}{c}{$\la=1$} & \mcol{1}{c|}{$\la=2$}\\
\hline
&&&\\[-0.3cm]
0  & $3.433(-2)$ & $9.869(-5)$ & $3.414(-2)$ \\
1  & $8.333(-4)$ & $9.869(-5)$ & $6.041(-4)$ \\
3  & $9.241(-5)$ & $9.869(-5)$ & $8.876(-5)$ \\
4  & $1.218(-7)$ & $2.279(-6)$ & $2.582(-6)$ \\
6  & $2.125(-8)$ & $4.987(-7)$ & $4.842(-7)$ \\
\hline
\end{tabular}
\end{center}
\end{table}

In Table 4 we present the absolute relative error in the computation of ${\cal W}^+_{\la,\nu}(x)$ using (\ref{e42}). The examples shown correspond to situations with $N=0, 1$ and 2. In the cases $N=1$ and 2 the value of $\Re h(u_N)<0$, so that $I_N(a,\la;x)$ is exponentially small and is neglected. The contributory series $I_j(a,\la;x)$ were calculated at optimal truncation, that is truncation at or near the term of least magnitude; this required the calculation of up to 30 coefficients $A_k$ and $B_k^{(1)}$.

Finally, to confirm the presence of the series $I_1(a,\la;x)$ in the $N=1$ case, we show in the upper half of Table 5 the values of  $\Delta{\cal W}\equiv {\cal W}^+_{\la,\nu}(x)-I_0(a,\la;x)$ for three values of $x$ compared with $I_1(a,\la;x)$, both series being optimally truncated. The first case $\la=3$, $a=0.20$ corresponds to $\Re h(u_1)<0$ whereas the second case $\la=4$, $a=0.20$ corresponds to $\Re h(u_1)>0$. In the lower half of the table we present similar calculations when $\la=6$, $a=0.20$ corresponding to $N=2$. In this case it is found that $\Re h(u_2)<0$ and so the exponentially small series $I_2(a,\la;x)$ is neglected.

\begin{table}[h]
\caption{\footnotesize{Values of the absolute relative error in ${\cal W}^+_{\la,\nu}(x)$ resulting from the expansion (\ref{e42}) for different truncation index $k$ when $x=20$.}}
\begin{center}
\begin{tabular}{|l||c|c|c|}
\hline
\mcol{1}{|c||}{$k$} & \mcol{1}{c|}{$\la=1,\ a=0.50$} & \mcol{1}{c|}{$\la=3,\ a=0.20$} & \mcol{1}{c|}{$\la=6,\ a=0.20$}\\
\mcol{1}{|c||}{} & \mcol{1}{c|}{$N=0$} & \mcol{1}{c|}{$N=1$} & \mcol{1}{c|}{$N=2$}\\
\hline
&&&\\[-0.3cm]
0  & $3.730(-3)$     &  $3.787(-3)$ & $1.550(-2)$\\
1  & $3.020(-6)$     &  $2.432(-4)$ & $3.381(-3)$\\
2  & $3.898(-6)$     &  $6.006(-5)$ & $2.963(-4)$\\
3  & $3.919(-7)$     &  $1.100(-5)$ & $1.256(-4)$\\
4  & $2.813(-8)$     &  $1.774(-6)$ & $8.332(-4)$\\
5  & $\,\,7.909(-10)$&  $1.786(-7)$ & $1.419(-5)$\\
\hline
\end{tabular}
\end{center}
\end{table}

\begin{table}[h]
\caption{\footnotesize{Values of $\Delta{\cal W}(x)\equiv {\cal W}^+_{\la,\nu}(x)-I_0(a,\la;x)$ compared with $I_1(a,\la;x)$ when $N=1$ and $N=2$ for different values of $x$.}}
\begin{center}
\begin{tabular}{|c|lll|}
\hline
&&&\\[-0.3cm]
\mcol{1}{|c|}{} & \mcol{3}{c|}{$\la=3,\ a=0.20\ \ (N=1)$} \\
\mcol{1}{|c|}{$x$} & \mcol{1}{c}{$\Delta{\cal W}(x)$} & \mcol{1}{c}{$I_1(a,\la;x)$} & \mcol{1}{c|}{${\cal W}^+_{\la,\nu}(x)$} \\
\hline
&&&\\[-0.3cm]
20  & $-1.58935(-2)$ & $-1.57281(-2)$ & $+7.070661(+5)$ \\ 
30  & $-1.48072(-2)$ & $-1.48186(-2)$ & $+1.986142(+9)$ \\
40  & $-8.74902(-3)$ & $-8.74792(-3)$ & $+5.920851(+12)$\\
\hline
&&&\\[-0.3cm]
\mcol{1}{|c|}{} & \mcol{3}{c|}{$\la=4,\ a=0.20\ \ (N=1)$}\\
\hline
&&&\\[-0.3cm]
20 & $-4.21656$ &  $-4.20876$  &       $+2.277758(+5)$\\
30 & $-3.00021(+1)$ & $-3.00057(+1)$ & $+3.823713(+8)$\\
40 & $-7.95934$(+2) & $-7.95905(+2)$ & $+6.805185(+11)$\\
\hline
&&&\\[-0.3cm]
\mcol{1}{|c|}{} & \mcol{3}{c|}{$\la=6,\ a=0.20\ \ (N=2)$} \\
\hline
&&&\\[-0.3cm]
20 & $+4.36797(+1)$ & $+4.31217(+1)$ & $+5.336787(+4)$\\
30 & $+1.45878(+4)$ & $+1.45867(+4)$ & $+4.687193(+7)$\\
40 & $-1.01722(+6)$ & $-1.01707(+6)$ & $+4.352648(+10)$\\
\hline
\end{tabular}
\end{center}
\end{table}

\end{document}